\documentclass{amsart}
\usepackage[arrow,matrix,curve,color]{xy}
\usepackage{color}
\usepackage[normalem]{ulem}
\usepackage{tikz}
\usepackage{xspace}
\definecolor{light-blue}{rgb}{0.8,0.85,1}
\definecolor{light-red}{rgb}{1,.4,.4}
\definecolor{purp}{rgb}{.7,.3,1}
\definecolor{yel}{rgb}{1,1,.5}
\definecolor{cy}{rgb}{0,1,1}
\usepackage{colortbl}
\usepackage{multirow}
\usepackage{array}
\makeatletter
\@namedef{subjclassname@2020}{%
  \textup{2020} Mathematics Subject Classification}
\makeatother

\theoremstyle{definition}

\newcommand{\bT}{\mathbb T}
\newcommand{\bR}{\mathbb R}
\newcommand{\bC}{\mathbb C}

\newcommand{\bH}{\mathbb H}

\newcommand{\bZ}{\mathbb Z}

\newcommand{\bP}{\mathbb P}

\newcommand{\cA}{\mathcal A}

\newcommand{\cH}{\mathcal H}

\newcommand{\cK}{\mathcal K}
\newcommand{\cL}{\mathcal L}

\newcommand{\Hom}{\operatorname{Hom}}

\newcommand{\ev}{\text{\textup{even}}}
\newcommand{\od}{\text{\textup{odd}}}
\newcommand{\dR}{\text{\textup{deR}}}

\newcommand{\co}{\colon\,}

\title{Twisted cohomology}
\author{Jonathan Rosenberg}
\address{Department of Mathematics\\
University of Maryland\\
College Park, MD 20742-4015, USA} 
\email{jmr@umd.edu}
\urladdr{http://www.math.umd.edu/~jmr/}

\begin{document}
\begin{abstract}
  We discuss twisted cohomology, not just for ordinary cohomology
  but also for $K$-theory and other exceptional cohomology theories,
  and discuss several of the applications of these in mathematical physics.
  Our list of applications is by no means exhaustive, but we are hoping
  that it is extensive enough to give the reader a feel for the
  possible applications of twisted theories in many different contexts.
  We also give many suggestions for further reading, but this subject
  has now expanded to the point where the bibliography is necessarily
  very incomplete.
\end{abstract}

\maketitle
\tableofcontents

\section{Twisted ordinary cohomology}
\label{sec:ord}
\subsection{Cohomology with local coefficients}
\label{sec:lc}
The story of twisted cohomology begins with the classical subject of
\emph{cohomology with local coefficients} \cite{MR9114}.
One way to motivate this
is to think about the Leray-Serre spectral sequence of a fibration
$F\to E\xrightarrow{p} B$, where $F$, $E$, and $B$ are nice topological spaces,
say CW-complexes.  When $B$ is simply connected,
the Leray-Serre spectral sequence
$H^p(B, H^q(F; G))\Rightarrow H^{p+q}(E; G)$ enables one in many
cases to compute
the cohomology of $E$ in terms of the cohomologies of $B$ and $F$.
However, when $B$ is path connected but not simply connected, the
fundamental group of $B$ can act nontrivially on the cohomology
$H^q(F; G)$ of the fiber, and $x\mapsto H^q(p^{-1}(x); G)$ can
become a nontrivial locally constant sheaf on the base $B$.
The cohomology of this sheaf is called cohomology with local coefficients.
In the case where $F$, $E$, and $B$ are the classifying spaces of
groups $K$, $L$, and $M$, and the fibration
$F\to E\xrightarrow{p} B$ comes from a group extension
$1\to K\to L\to M\to 1$, the Leray-Serre spectral sequence specializes
to the Hochschild-Serre spectral sequence of the group extension
$H^p(M, H^q(K; G))\Rightarrow H^{p+q}(L; G)$, where one has to take
into account the action of the quotient group $M$ on the cohomology
of the normal subgroup $K$.

Another important place where cohomology with local coefficients
is needed is in expressing Poincar{\'e} duality for non-orientable
manifolds.  Suppose $M^n$ is a connected compact manifold.
If $M$ is orientable,
Poincar{\'e} duality says that cap product with the fundamental
homology class $[M]\in H_n(M;\bZ)$ (for some choice
of orientation) gives an isomorphism from
$H^k(M;\bZ)$ to $H_{n-k}(M;\bZ)$.  When $M$ is
not orientable, this can no longer work since there is no
fundamental homology class $[M]\in H_n(M;\bZ)$; in fact $H_n(M;\bZ)=0$.
In this case, the first Stiefel-Whitney class
$w_1\in H^1(M;\bZ/2)$ is non-zero. We can repair Poincar{\'e} duality
by observing that there is a twisted fundamental class $[M]$ in $H_n(M)$
\emph{with local coefficients} given by $w_1$, and that cap
product with $[M]$ implements an isomorphism from
$H^k(M;\underline{\bZ})$ to $H_{n-k}(M;\bZ)$, where $\underline{\bZ}$
denotes the locally constant sheaf locally isomorphic to $\bZ$
on which $\pi_1(M)$ acts via $w_1(M)$.  (Think of $w_1$ as taking
values in $\{\pm 1\}$ and let it act on $\bZ$ by multiplication.)
This was first observed by Steenrod in \cite{MR9114} and has played
a big role in manifold theory ever since.

\subsection{Twisted de Rham cohomology}
\label{sec:twisteddR}
Cohomology with local coefficients
is the first example of twisted
cohomology.  But there are other twists of ordinary cohomology, even
of de Rham cohomology on manifolds (which is naturally isomorphic
to ordinary cohomology
with coefficients in $\bR$ or $\bC$). The following construction
seems to have originally been ``folklore,'' but a hint of it appears
in \cite{MR2091892}, with a longer discussion in \cite[\S9.3]{MR1911247},
and a full exposition in \cite{Twatson}.

Let $M$ be a manifold and let
$\alpha$ be a closed $k$-form on $M$, where $k$ is \emph{odd}.
Then $d'_\alpha=d+\alpha\wedge$ maps even forms to odd forms and odd forms to
even forms. (However, it increases degree by $1$ only
when $\alpha$ is a $1$-form, i.e., $k=1$.)
And $d'_\alpha$ is a differential, i.e.,
\[
(d'_\alpha)^2\omega = (d+\alpha\wedge\,)(d\omega+\alpha\wedge\omega)
= d^2\omega - \alpha\wedge d\omega + \alpha\wedge d\omega
+ (\alpha\wedge\alpha)\wedge\omega = 0
\]
for any form $\omega$.  Thus the cohomology $H^\bullet_\alpha(M)$
of the $2$-periodic complex
\[
\cdots \xrightarrow{d'_\alpha} \Omega^\ev(M) \xrightarrow{d'_\alpha} \Omega^\od(M)
\xrightarrow{d'_\alpha} \Omega^\ev(M) \xrightarrow{d'_\alpha} 
\cdots
\]
can be called the \emph{twisted cohomology} of $M$ with twist $\alpha$.
A simple calculation \cite[Theorem 5]{Twatson} shows that
$d'_\alpha$ and $d'_{\alpha'}$ are conjugate under $e^{-\eta}$ if
$\alpha'-\alpha=d\eta$, so $H^\bullet_\alpha(M)$ only depends on
the de Rham class of $\alpha$.
Note that (unless $\alpha$ is a $1$-form) we will only think of it
as being $\bZ/2$-graded rather than $\bZ$-graded, but when
$\alpha=0$, it coincides with the usual de Rham cohomology of $M$,
$\bZ/2$-graded via even and odd degree forms.  Also note (see
\cite[\S3.3]{Twatson}) that $d'_\alpha$ is the total differential of
a double complex, where the horizontal differential is the usual
differential $d$ and the vertical differential is $\alpha\wedge\,$.
In the spectral sequence of the double complex, the first differential
$d_1$ is $d$ and so $E_2$ is just
the usual de Rham cohomology $H^\bullet_\dR(M)$, and the
next differential $d_k$
is cup product with the de Rham class $[\alpha]$ of $\alpha$.  So
the next stage $E_{k+1}$ in the spectral sequence for
$H^\bullet_\alpha(M)$ is the cohomology of $[\alpha]\cup \underline{\phantom{x}}$
acting on $H^\bullet_\dR(M)$.  However, in general there can be still
higher differentials.  Atiyah and Segal \cite[Proposition 6.1]{MR2307274}
have shown that these higher differentials all come from iterated
Massey products with $\alpha$.

There is one case where twisted cohomology reduces to cohomology with
local coefficients.  Assume $k=1$ so that $\alpha$ is a $1$-form.
Since $H^1_{\dR}(M)\cong H^1(M;\bR)\cong \Hom(\pi_1(M),\bR)$,
$[\alpha]$ can be identified with a group homomorphism $\pi_1(M)\to \bR$.
There is an associated locally constant sheaf $\cL$ and one can easily
check that $H^\bullet_\alpha(M)\cong H^\bullet(M, \cL)$, which is cohomology with
local coefficients \cite[\S3.1]{Twatson}.

The kind of twisted cohomology we have just been discussing is based
on de Rham cohomology, so it ignores torsion phenomena.  In
\cite{MR2313935} is a proposal for a kind of integral twisted
cohomology based on homotopy theory.  Namely, given a cohomology class
$\alpha\in H^k(X;\bZ)$ (the main case of interest is when $k=3$), one
can identify $\alpha$ with the classifying invariant for a principal
fibration $K(\bZ,k-1) \to E \to X$.  The Serre spectral sequence for
this fibration will have $E_2^{p,0} = H^p(X;\bZ)$ and first
differential $d_k$ sending the generator $\iota$ of
$H^{k-1}(K(\bZ,k-1); \bZ)\cong \bZ$ to $\alpha\in H^k(X;\bZ)$, and
sending $x\otimes \iota$ to $\pm\alpha\cup x$, just as for
the spectral sequence in twisted de Rham cohomology. It is proposed in
\cite{MR2313935} to take the $\alpha$-twisted cohomology of $X$ to be
the cohomology of $E$, though this has the drawbacks that if
$\alpha=0$, it computes $H^\bullet(X\times K(\bZ,k-1))$ and not
$H^\bullet(X)$, and even if $X$
is finite-dimensional, this will usually have cohomology in arbitrarily
large degree.  So it seems a more reasonable definition would be the
image $E_\infty^{p,0}$ of the edge homomorphism $H^p(X;\bZ)\to H^p(E;\bZ)$,
which agrees with $H^p(X;\bZ)$ when $\alpha=0$ and always vanishes beyond
the dimension of $X$.

\section{Twisted generalized cohomology theories}
\label{sec:gentheories}

\subsection{``Classical'' twisted $K$-theory}
\label{sec:twistedK}

When it comes to recent impact on mathematical physics, one can easily
argue that twists of $K$-theory and other generalized cohomology theories
(not satisfying the Eilenberg-Mac Lane dimension axiom) have played
a much bigger role than twisted ordinary cohomology.  We shall begin with
$K$-theory, the most familiar and important generalized cohomology theory,
and then move on later to other theories.

The story of twisted $K$-theory begins with the work of Karoubi and Donovan
\cite{MR256386,MR282363}. They actually worked with $\bZ/2$-graded algebras,
but for simplicity we'll ignore this here. Then
they observed that if one looks at Azumaya
algebras over a compact Hausdorff space $X$ (algebras of the form
$A=\Gamma(X,\cA)$, where $\cA$ is a locally trivial bundle over $X$ with fibers
of the form $M_n(D)$, the $n\times n$ matrices
over a finite-dimensional division algebra $D$ over $\bR$ or $\bC$,
so that $D\cong \bR$, $\bC$ or $\bH$), the $K$-theory of the algebra
$A$ looks like a version of 
$KO^\bullet(X)$, $KU^\bullet(X)$, or $KSp^\bullet(X)$ (depending on which division
algebra one is using) with
local coefficients, where the local coefficients come from the twisting
of the bundle $\cA$. In, say, the complex case where $\cA$ has
fibers $M_n(\bC)$, $\cA$ is classified by the homotopy class of a map from
$X$ to the classifying space $BPGL_n(\bC)\simeq BPU(n)$,
and because of the Morita
invariance of $K$-theory, the $K$-theory of $A$
only depends on the Brauer group class of the Azumaya algebra, which
is a torsion class in $H^3(X;\bZ)$.  (In fact, by a theorem of Serre,
every torsion class in $H^3(X;\bZ)$ arises this way.)

The theory of \emph{continuous-trace algebras}, for which a good reference
is \cite{MR1634408}, provides a substitute for Azumaya algebras when the
characteristic class in $H^3(X;\bZ)$, now called the \emph{Dixmier-Douady
  class}, is no longer torsion.  A ``stable'' continuous-trace algebra, say
over a connected second-countable locally compact Hausdorff space $X$, is a
$C^*$-algebra of the form $A=\Gamma_0(X,\cA)$, where $\cA$ is a locally trivial
bundle over $X$ with fibers isomorphic to $\cK$, the compact operators on
a separable infinite-dimensional
Hilbert space $\cH$ (over $\bR$, $\bC$ or $\bH$), and structure group
the $C^*$-algebra automorphisms of $\cK$.  Such an algebra
is classified by the homotopy class of a map from
$X$ to the classifying space $BPU(\cH)$.  Since the unitary group $U(\cH)$
is contractible (in the natural topology, the strong operator topology),
this classifying space has the homotopy type of the Eilenberg-Mac Lane space
$K(\bZ, 3)$ in the complex case (since the center of $U(\cH)$ is the
circle group $\bT$ and thus $PU(\cH)$ is a $K(\bZ,2)$), or $K(\bZ/2, 2)$
in the real or quaternionic cases (since then the center of $U(\cH)$
consists just of $\{\pm1\}$, and thus $PU(\cH)$ is a $K(\bZ/2,1)$).
For details on all of this see \cite{MR1018964}.  Thus given
$\delta\in H^3(X;\bZ)$, one has twisted $K$-groups $K^\bullet_\delta(X)$
in the complex case, and given $\delta\in H^2(X;\bZ/2)$, one has
twisted $K$-groups $KO^\bullet_\delta(X)$ and $KSp^\bullet_\delta(X)$
in the real and quaternionic cases, obtained by taking the
topological $K$-theory of the appropriate continuous-trace algebras.

There is a slight technical complication in that to get a functorial map
$K^\bullet_{\delta'}(Y)\to K^\bullet_\delta(X)$, one needs a map 
$f\co X\to Y$ that induces a $*$-homomorphism of the appropriate
$C^*$-algebras, and thus a map of bundles.  This involves a choice
of \v Cech cocycles representing $\delta$ and $\delta'$ such that $f$ carries
one cocycle to the other.  More formal definitions of twisted $K$-theory
and discussions of this and other fine points may be found in
\cite{MR2172633,MR2307274} and \cite{MR2513335}.

As one would expect, there is a compatibility between twisted $K$-theory
and twisted ordinary cohomology, given by the Chern character.  If $M$
is a compact manifold (the compactness is not so essential, but is
assumed here to avoid complications if $M$ is not of finite type)
and $\delta\in H^3(M;\bZ)$, then
the (twisted) Chern character gives a natural transformation from
$K^\bullet_\delta(M)$ to $H^\bullet_\delta(M)$ which is an isomorphism after
tensoring the $K$-theory with $\bR$.  (Here one should take a closed
$3$-form representing the integral cohomology class $\delta$.)
For references, see \cite{MR1911247}, \cite{MR1977885}, and
\cite{MR2307274}.

One immediate use for twisted $K$-theory, pointed out in the original
work by Karoubi and Donovan, involves the Thom isomorphism.  Suppose
$X$ is a compact Hausdorff space and $V$ is a real vector bundle of
rank $k$ over $X$.  If $V$ is orientable, a choice of orientation on $V$
gives rise to a Thom isomorphism
$\tau\co H^\bullet(X)\xrightarrow{\cong} H^{k+\bullet}_c(V)$,
where the cohomology on the right is taken with compact support.
(Alternatively, one could take the relative cohomology of the disk bundle
of $V$ relative to its boundary the sphere bundle, or take the
reduced cohomology of the Thom space.)  When $V$ is not orientable,
one gets a Thom isomorphism only if the cohomology of $V$ (with compact
support) is replaced by twisted cohomology with twist $w_1$, $w_1$ the
first Stiefel-Whitney class of the bundle.

This suggests a similar picture for $K$-theory.  (Here we only deal with
complex $K$-theory for simplicity, but the same would work for $KO$
with spin$^c$ replaced by spin.)  Suppose $V$ is equipped with a
spin$^c$ structure.  Then again one gets a Thom isomorphism
$\tau\co K^\bullet(X)\xrightarrow{\cong} K^{k+\bullet}(V)$,
where the $K$-theory is taken with compact support. (A proof of
this when $V$ is a complex vector bundle is in
\cite[Theorem 4.7]{MR228000}.)  If $V$ does
\emph{not} admit a spin$^c$ structure, then this fails, but
it can be repaired by replacing $K^{k+\bullet}(V)$ with
twisted $K$-theory.  For example, if $V$ is oriented but not spin$^c$,
one should use
$K^{k+\bullet}_{W_3(V)}(V)$, the twisted $K$-theory for the obstruction
$W_3(V)\in \text{Tors}\,H^3(V;\bZ)$ to a spin$^c$ structure on $V$
(see \cite[\S6]{MR282363} or \cite[\S4]{MR2513335}).

In general it is not so obvious how to compute twisted $K$-theory,
but there are two useful techniques in this direction.  One is the
(twisted) Atiyah-Hirzebruch spectral sequence
$H^p(X,K^q)\Rightarrow K^{p+q}_\delta(X)$.  This has the same $E_2$
term as the Atiyah-Hirzebruch spectral sequence for untwisted $K$-theory,
but the difference is that the twist $\delta$ appears in the first
nonzero differential $d_3$, which is now
$\text{Sq}^3+\delta\cup\underline{\phantom{x}}$
\cite{MR679694,MR2172633}.  A second quite different technique
is to use a theorem of Khorami \cite{MR2832567}, which shows that
for $\delta\in H^3(X;\bZ)$ and $P$ the associated $PU(\cH)$-bundle
over $X$, one has
$K_\bullet^\delta(X)\cong K_\bullet(P)\otimes_{K_0(\bC\bP^\infty)}\bZ$,
where the map $K_0(\bC\bP^\infty)\to \bZ$ is rather subtle and
is not flat as a map of rings. (The ring structure on $K_0(\bC\bP^\infty)$
comes from the Pontryagin product.)  Note that this result
is formulated in $K$-twisted homology, but one can dualize
via the universal coefficient theorem (which applies to
continuous-trace algebras \cite{MR731763}) to compute twisted $K$-cohomology.

\subsection{``Exotic'' twisted $K$-theory}
\label{sec:exotictwistedK}

What we have described so far is what one could call
``classical'' twisted $K$-theory, known since the work of Karoubi and Donovan
in 1969--1970.  At this point quite a number of generalizations are known.
First of all, there is an equivariant theory for spaces equipped
with an action of a compact Lie group $G$.  This is discussed, for example,
in \cite[\S6]{MR2172633}or \cite[\S6]{MR2513335}.  In the complex
(ungraded) case, the twisting is by a class in equivariant cohomology
$H^3_G(X;\bZ)$.

A second variant is \emph{twisted $KR$-theory}, developed in
the work of Moutuou \cite{Moutuou,MR3177819} and in
\cite{MR3267662}.  To explain this let's start with a review of
Atiyah's $KR$-theory \cite{MR206940}.  This is a cohomology
theory based on vector bundles
on ``Real spaces,'' (locally compact) spaces $X$ equipped with
an involutive homeomorphism $\iota$.  The vector bundles are
required to have a conjugate-linear involution compatible with the
involution on the base space $X$.  The motivating example of a Real space is
the set of complex points of an algebraic variety defined over $\bR$,
with $\iota$ given by the action of $\text{Gal}(\bC/\bR)$.
Like $KO$, $KR$-theory is $8$-periodic, and when $\iota=\text{id}_X$,
$KR^\bullet(X,\iota)=KO^\bullet(X)$.  On the other hand, when $X=Y\sqcup Y$
and $\iota$ interchanges the two copies of $Y$, a Real vector bundle on
$X$ is uniquely determined by a complex vector bundle on $Y$, so that
$KR^\bullet(X)\cong K^\bullet(Y)$.

However, as we have already seen in the discussion of twisted real
$K$-theory, at some level the division algebras $\bR$ and $\bH$ have
to be treated on an equal footing.  This insight is reinforced by the
fact, due to Donald Anderson, that $KO^\bullet$ and $KSp^\bullet$ are naturally
duals of each other, whereas $K^\bullet$ (complex $K$-theory) is self-dual.
This motivates the consideration \cite{MR3267662} of a twisted version
of $KR$ theory, in which for $(X,\iota)$ a Real space, each connected
component of the fixed set $X^\iota$ of $\iota$
is assigned a sign $\pm$.  When all the signs are $+$,
one recovers the usual $KR^\bullet$, but for each $-$ sign, the restriction
of $KR_\pm^\bullet (X)$ to the associated component of $X^\iota$ is
$KSp^\bullet$ rather than $KO^\bullet$.  Still another version of
twisted $KR$-theory, with some applications to Grothendieck-Witt
groups of real algebraic varieties, may be found in
\cite{MR3614970}.

There is another way of twisting $K$-theory which, even in the case
of non-equivariant complex $K$-theory, produces something new.
This approach was first introduced in \cite{MR2681757} with special
attention to twisted $K$-theory, and was also studied in
\cite{MR4163521}.  In the next subsection we will see how this generalizes
to other cohomology theories.

For simplicity we just deal with complex $K$-theory.  As is well known,
the classifying space for stable complex vector bundles, which is
also the classifying space for $K$-theory, is $\bZ\times BU$, where
the $\bZ$ factor is used to record the rank of the bundle and $BU$ is
the classifying space of the infinite unitary group.  Twisting by a class
in $H^3$ comes, as we have seen before, from the identification
of $BPU(\cH)$, the classifying space for the projective unitary
group of a separable infinite-dimensional Hilbert space,
with a $K(\bZ,3)$ space, and the fact that one can replace $U$ by
unitaries of the form $1+(\text{compact})$, which are preserved under
conjugation by $PU(\cH)$.  However, $PU(\cH)$ does not account for the
full set of homotopy automorphisms of $K$-theory.  To study these
one should build a topological monoid $GL_1(K)$, which in other language
was computed in \cite{MR494076} --- see \cite[Chapter 6]{MR505692}
and \cite{MR3192610} for an explanation.  It turns out that
\[
GL_1(K)\simeq K(\bZ/2, 0) \times K(\bZ, 2) \times BSU^\otimes
\]
and thus
\[
BGL_1(K)\simeq K(\bZ/2, 1) \times K(\bZ, 3) \times BBSU^\otimes.
\]
Here the first two factors are ``classical'': the $K(\bZ/2, 0)$
comes from complex conjugation on $U$ and the $K(\bZ, 2)$ comes from
$PU(\cH)$.   But the last factor is more mysterious, and leads to
new twists of $K$-theory.  $BSU^\otimes$ is the
infinite loop space classifying virtual complex vector bundles
of rank and determinant one, equipped with the tensor product structure.
After localization at any prime $p$, it is equivalent to the usual
$BSU$, so it has lots of nonzero homotopy groups.  Then a homotopically
nontrivial map $X\to BBSU^\otimes$ can be used to construct an
exotic twist of $K$-theory, different from the twists that come
from elements of $H^3(X; \bZ)$.

\subsection{Twisted theories in general}
\label{sec:gentheory}

Now that we have examined a large number of examples of twisted cohomology
or cohomology with local coefficients, it is time to discuss a general
framework for such theories, and to explain how this framework applies
to cohomology theories other than ordinary cohomology and $K$-theory.
This theory was developed by Ando, Blumberg, and Gepner
\emph{et al.} in \cite{MR3252967,MR3286898,MR3890766}.
Alternative approaches due to May and Sigurdsson and to
Hebestreit, Sagave, and Schlichtkrull may be found in
\cite{MR2271789,MR4079640,MR4163521}.  To avoid very complicated homotopy
theoretic problems, the treatment here will be a bit sloppy, but should
convey the general idea.  For details of how to make everything rigorous,
please see the original references.

Let's start with a homotopy-theoretic interpretation of cohomology with
local coefficients.  Recall that if $G$ is an abelian group, then
$H^k(X;G)$ can be identified with the set of homotopy classes of maps
$X\to K(G,k)$ into the Eilenberg-Mac Lane space $K(G,k)$.  For cohomology
with local coefficients $H^k(X, \underline{G})$, $\underline{G}$ a local
coefficient system coming from an action of $\pi_1(X)$ on $G$,
one can similarly identify the local cohomology with 
$\pi_1(X)$-\emph{equivariant} homotopy classes of equivariant maps from
the universal cover $\widetilde X$ of $X$ to $K(G,k)$.  This approach
won't generalize to more general twisted theories, but we can instead
create a fiber bundle
\[
E=\widetilde X\times_{\pi_1(X)} K(G,k) \to \widetilde X\times_{\pi_1(X)} * = X
\]
over $X$
with fibers that are copies of $K(G,k)$.  Then $H^k(X, \underline{G})$
can be identified with homotopy classes of sections of this bundle.

A similar approach will work with more general cohomology theories.
A generalized cohomology theory $E$ is always represented in the
homotopy category by a \emph{spectrum}, which is a sequence
$\{E_k\}_{k\ge 0}$ of (based) spaces with connecting maps
$f_k\co \Sigma E_k\to E_{k+1}$.  Then we define
\[
E^n(X) = \varinjlim_k [\Sigma^{k-n} X, E_k], 
\]
where $[\underline{\phantom{x}}, \underline{\phantom{x}}]$ denotes
homotopy classes of (based) maps and we map
$f\co \Sigma^{k-n} X\to E_k$ to $\Sigma f\co \Sigma^{k-n+1} X\to \Sigma E_k$
and then compose with $f_k$ to get a map $\Sigma^{k-n} X \to E_{k+1}$.
The definition
is rigged so that $E^{n+1}(\Sigma X)\cong E^n(X)$, which is one of the
key properties of a cohomology theory.

For present purposes we will mostly be interested in generalized
cohomology theories $E$ that come with a suitable notion of cup
products.  In this case, we need some extra structure on $E$,
summarized by saying that it is a \emph{ring spectrum}, i.e., is
equipped with a suitable multiplication map $E\wedge E\xrightarrow{\mu_E} E$
that is homotopy associative.  There
are various axioms to satisfy, but these hold in most of the
standard examples, such as ordinary cohomology with coefficients
in an abelian group $G$ (where $E$ is
what's called an \emph{Eilenberg-Mac Lane spectrum} $HG$, where
$E_k = K(k,G)$ and $\Sigma K(k,G)\to  K(k+1,G)$ is adjoint to the
equivalence $K(k,G)\xrightarrow{\simeq} \Omega K(k+1,G)$),
real or complex $K$-theory, and oriented, unoriented, spin, and
complex cobordism.  Given a ring spectrum $E$, we also have the notion
of an $E$-\emph{module spectrum}, which is a spectrum $F$ with a multiplication
map $\mu_F\co E\wedge F\to F$ such that
\[
\xymatrix@C+2pc{E\wedge E\wedge F
  \ar[r]^(.55){{\mu_E}\wedge \text{id}_F}\ar[d]^{\text{id}_E\wedge \mu_F}
  & E\wedge F\ar[d]^{\mu_F}\\
  E\wedge F\ar[r]^(.55){\mu_F} & F}
\]
is homotopy commutative.

Twisted $E$-cohomology of $X$ comes from taking a bundle of
$E$-module spectra
over $X$, with fibers that are equivalent to $E$, and taking
homotopy classes of sections.  In the case of a trivial bundle
$X\times E\xrightarrow{\text{pr}_1} X$, a section is just a 
map $f \co X\to E$, so we recover the homotopy group $[X, E]$,
which is $E^0(X)$, and similarly after suspending for other values of $n$.
The sort of bundle needed is what is called a ``parameterized spectrum''
in the literature.  The set of these bundles is parameterized by
homotopy classes of maps $X\to BGL_1(E)$, so this definition recovers
our original definition of cohomology with local coefficients or
twisted $K$-theory.  An essentially equivalent definition involves
giving a more general notion of Thom space (a ``Thom spectrum'') and
taking its homotopy groups.

\section{Some physical applications}
\label{sec:appl}

Let's now outline a few of the major applications of twisted cohomology
in mathematical physics.

\subsection{D-brane charges in Type II string theory}
\label{sec:Dbranes}

Our first example comes from type II string theory.  In this theory,
spacetime is a $10$-dimensional manifold $X$ (usually spin)
and it contains submanifolds
$N\subseteq X$ called \emph{D-branes}, of even dimension in type IIB and
of odd dimension in type IIA, on which ``open'' strings can begin and
end.  (The letter D stands for ``Dirichlet,'' from Dirichlet boundary
conditions.) These D-branes carry \emph{Chan-Paton bundles} and
Ramond-Ramond charges.  It is now generally believed that D-brane charges
should take values in $K$-theory (of even degree in type IIB,
of odd degree in type IIA) rather than in ordinary cohomology
\cite{MR1606278,MR1621204,MR1674715}. However, this statement is
really only true when the $B$-field of string theory is topologically
trivial (i.e., given by a closed $2$-form).  If this is not the
case, then there is an ``$H$-flux'' $H\in H^3(X;\bZ)$, and D-brane
charges should take their values in $H$-twisted $K$-theory
\cite{MR1674715,MR1807598,MR1756434}.

Let's give a quick review of the physical argument for the
need for twisted $K$-theory.  The key fact is that the Freed-Witten
anomaly formula on open strings \cite{MR1797580} forces, for $N$
a D-brane in spacetime $X$, $W_3(N)$ to agree with the pullback
of $H$ from $X$ to $N$. To phrase things another way (see
\cite[\S4]{MR2079376} or \cite{MR1807598}), neither the gauge field
associated to the Chan-Paton bundle nor the $B$-field as a $2$-form are
completely well-defined on the D-brane, but the indeterminacies
have to cancel so that the path integral is well defined.
When $H=0$, that means that $N$ is a spin$^c$
manifold; in fact, in classical IIB string theory, one often thinks
of the D-branes as being complex submanifolds of $X$ (which are of
course automatically spin$^c$).  But in general, the D-branes are not
spin$^c$, so the Poincar\'e dual of the Chan-Paton class in $K$-homology
is not well-defined.  However, the fact that $N$ is ``twisted spin$^c$''
means that one does have a Poincar\'e dual in $H$-twisted $K$-homology,
which one can push forward to a class in $K^H_\bullet(X)$ or to
its dual $K_H^\bullet(X)$.

A slightly different use of twisted $K$-theory in string theory
comes from an attempt to relate type IIB string theory to the
$E_8$ gauge theory in eleven dimensions that underlies M-theory
\cite{MR2062363}.  In this picture the twisting of $K$-theory
in $H^3(X;\bZ)$ can be identified with the homotopy class of a
map $X\to E_8$ (since $E_8$ is a good approximation to a
$K(\bZ,3)$ space) and since $[X,E_8]=[X,\Omega BE_8]\cong [\Sigma X, BE_8]$
we obtain an $E_8$ bundle over a circle bundle $Y$ over $X$.

\subsection{Topological T-duality}
\label{sec:toopTdual}

The concept discussed above in subsection \ref{sec:Dbranes}, that
D-brane charges naturally live in twisted $K$-theory, has an interesting
consequence which can be used to test various physical hypotheses.
Namely, physicists predict many ``dualities'' in which seemingly
different mathematical models give rise to the same physical
phenomena.  An example of such dualities is T-duality
(or target-space duality or torus duality) in type II string theory,
in which two different type II string theories on possibly different
spacetimes $X$ and $X'$, obtained by replacing certain circles
in $X$ by their duals, are actually physically equivalent.  It was
first observed in \cite{MR2080959,MR2116165} that under T-duality,
in general both the topology of $X$ and the $H$-flux in $H^3(X;\bZ)$
have to change.  Twisted $K$-theory imposes a constraint on this
change of topology; since the D-brane charge group is given by
twisted $K$-theory, under a T-duality between $(X,H)$ and $(X',H')$,
one needs to have $K^\bullet_H(X)\cong K^{\bullet+k}_{H'}(X')$, where
$k$ is the number of T-dualities involved.

The Bouwknegt-Evslin-Mathai observation led to a fairly large
industry to studying ``topological T-duality,'' trying to understand
the algebraic topology underlying the physics of T-duality.
This began with twisted cohomology \cite{MR2062361}
and moved on to a full understanding of how twisted $K$-theory
is related to the $H$-flux on a principal torus bundle
\cite{MR2116734,MR2130624,MR2287642,MR2222224,MR2482327,MR2491618,MR2560910}.
Interesting phenomena that occur in the case of higher-dimensional
torus bundles are that sometimes T-duals don't exist, are non-unique,
or are noncommutative.  The interested reader will find a very
rich literature on these topics, of which we have only scratched
the surface.

Still another variant of this situation comes from studying
topological T-duality for \emph{non-principal} circle or
torus bundles. For this
variant of the theory see \cite{MR2985336,MR3282986,MR3342758,MR3285613}.

Or one can replace bundles by ``bundles with degeneration'' or stacks.
This generalization of topological T-duality is studied in
\cite{MR2246781,MR2369414,MR2797285,MR2989461,MR3952356}.

Another generalization of topological T-duality is what has been
called ``spherical T-duality,'' in which circle bundles are replaced
by pricipal $SU(2)$-bundles.  (Recall that $SU(2)$ is topologically
just $S^3$.)  Such bundles are classified by homotopy classes of maps
into $BSU(2)\cong \bH\bP^\infty$.  In the simplest cases where one
considers $SU(2)$-bundles over a closed oriented $4$-manifold $M$,
these bundles are parametrized by $H^4(M;\bZ)\cong \bZ$ and the total
space of each bundle is a closed oriented $7$-manifold $P$.  In
\cite{MR3339166,MR3328056,MR3850273} it is shown that one gets
isomorphisms of twisted cohomology between $(P,H)$, $H\in H^7(P;\bZ)$,
and a dual bundle $(\widehat P,\widehat H)$,
$\widehat H\in H^7(\widehat P;\bZ)$.  There is also an isomorphism of
twisted $K$-theory with twisting $H$ in the more general sense of
\cite{MR2681757}.  When the base of the bundle has dimension bigger
than $4$, or one considers non-principal $SU(2)$-bundles,
then the situation is more complicated and spherical T-duals
may not exist or may be non-unique.

Finally, there are some extensions of topological T-duality that
interact with other subjects such as Fourier-Mukai and Langlands
duality.  For a few results in this direction see
\cite{MR2309993,MR3285609,MR3361543,MR3764065}.

\subsection{The Verlinde formula and the work of Freed-Hopkins-Teleman}
\label{sec:FHT}

One of the most striking applications of twisted (equivariant) $K$-theory
to mathematical physics is the work of Freed, Hopkins, and Teleman
\cite{MR2365650,MR2860342,MR3037783,MR2831111} on
using twisted $K$-theory to explain and understand the ``Verlinde formula''
\cite{MR954762}.  The latter, which on the face of it has nothing
to do with twisted $K$-theory, computes the structure constants for the
``fusion rule'' for multiplication of the primary fields in the
conformal field theory on $\bC\bP^1$ with three punctures coming from
a positive energy (projective) representation of the loop group $LG$ of
a compact Lie group $G$.  Roughly speaking, their result is that a
positive energy representation determines a twist $\tau$ in
$H^3_G(G;\bZ)$ (where $G$ acts on itself by conjugation), the sum of
the level of the representation and the dual Coxeter number, and the
fusion ring can be identified naturally with $^\tau K^\bullet_G(G)$,
the $\tau$-twisted equivariant $K$-theory, with ring operation coming
from Pontryagin product.   The actual details are more complicated and
include many auxiliary results of independent interest.
 
\subsection{Brane charges in the WZW model} 
\label{sec:WZWbranes}

Related to the Freed-Hopkins-Teleman theorem is the computation of the
non-equivariant twisted $K$-theory of a compact Lie group $G$, which
should classify the D-branes in the WZW theory of a given level.
A first attack on the problem of computing these twisted $K$-groups
was made in \cite{MR1834409,MR1877986,MR1960468,MR2045888,MR2080884},
followed by a more comprehensive
approach by Braun \cite{MR2061550}, still partially conjectural since
it relied on a conjecture about the commutative algebra of Verlinde
rings (which is known in some cases but may still be open in others).
Finally, using rather complicated techniques from algebraic topology,
Douglas proved Braun's result unconditionally in \cite{MR2263220},
though in a different form.  The result of \cite{MR2061550}
is that for $G$ a simply
connected simple compact Lie group of rank $n$ and twisting $h>0$ in
$H^3(G;\bZ)\cong \bZ$, $K^h_\bullet(X)$ is (as a ring) the tensor
product of an exterior algebra over $\bZ$ on $n-1$ odd-degree generators with a
finite cyclic group of order $c(G, h)=h/\text{gcd}(h,y(G))$, where
the constant
$y(G)$ can be described easily in terms of $n$ for the classical
groups (for example it is $\text{lcm}(1,2,\cdots,n)$ for $A_n=SU(n+1)$)
and where it takes the value $60$ for $G_2$, $27720$ for $F_4$
and $E_6$, $\text{lcm}(1,2,\cdots,17)$ for $E_7$, and
$\text{lcm}(1,2,\cdots,29)$ for $E_8$.
Other proofs of many of these results were given in \cite{MR3764065}
and \cite{MR4071368}.  Some of these calculations illuminate
certain conjectured physics dualities such as ``level-rank duality''
or ``strange duality.''

\subsection{Orientifold brane charges}
\label{sec:orient}

\emph{Orientifold string theories} are a variant of string theory
where one assumes that spacetime $X$ is equipped with an involution $\iota$,
and the map $\Sigma\to X$ of a string worldsheet into $X$ is assumed
to be equivariant for the \emph{worldsheet parity operator} $\Omega$
on $\Sigma$ and $\iota$ on $X$.  Such theories include type I string
theory, for example, obtained from type IIB string theory on $X$
by taking $\iota=\text{id}_X$.  Components of the fixed set $X^\iota$
are usually called O-planes (``O'' for orientifold).  As explained
in \cite[\S5.2]{MR1674715} and in \cite{MR1736794,MR1777343},
D-brane charges in orientifold theories should take their values in
$KR$-theory.  However, the $KR$-theory can be twisted because the
gauge field on an O-plane can be of either orthogonal or symplectic
type, leading to a need for sign-twisted $KR$-theory as discussed above in
subsection \ref{sec:exotictwistedK}.  Examples of calculations
with this theory in physically relevant situations were given in
\cite{MR3267662}.  Then in \cite{MR3316647}, a complete analysis was made of
all the possible twisted orientifold string theories on an elliptic curve
(with holomorphic or antiholomorphic involution), and the twisted
$KR$-theory was computed for each one.  Isomorphisms between the various
groups (further studied in \cite{MR3305978}) agree with physics-based
predictions of dualities between the various theories.  In these
calculations, one special case of twisted $KR$-theory that occurs is
actually twisted $KO$-theory in the sense of Karoubi and Donovan.
This case comes up
in the treatment of Witten's theory of ``type I string theory without
vector structure,'' introduced in \cite{MR1615617}.

\subsection{Applications in condensed matter theory}
\label{sec:condensed}

Twisted $K$-theory appears not only in string theory but also in
condensed matter physics.  The story begins with the ``tenfold way''
approach to classification of topological states of matter via $K$-theory
\cite{Kitaev,MR4575939}.  Most of the physical properties of solids
come from the quantum mechanics of the electrons that are not tightly bound
(those with energies
above the Fermi level).  Depending on the dimension and the
symmetry class (whether or not one has time-reversal symmetry,
parity symmetry, etc.),
the wave functions of these electrons can be viewed as sections of
a vector bundle with certain symmetry properties, and thus live
in a certain $K$-group that classifies the states.  It was shown in
\cite{MR3119923} that a more refined analysis (for example, taking
into account the various crystallographic groups), leads sometimes to
twisted $K$-theory, which can help explain phenomena such as topological
insulators.  By now there is quite a large mathematical physics literature
involving twisted cohomology or twisted $K$-theory in the study of
topological states of matter.  A few representative papers
are \cite{Barkeshli:2014cna,Chiu:2015mfr,MR3604557,MR3872629,MR4483417,MR4311720,MR4268163,MR4340936,MR4403305,MR4568912,MR4612787,MR4647667}.

\subsection{Applications of other twisted cohomology theories}
\label{sec:other}

Up to this point, most of the physics applications we have discussed
have concerned twisted $K$-theory of one form or another.  In this
last subsection we mention some applications of twists of other
cohomology theories.  One example that seems to have physical applications
is twisted Morava $K$-theory $K(2)$ at the prime $2$, or $2$-local
$\text{\textit{tmf}}$ (topological modular forms), which is closely related.
At the prime $2$, $K(2)$ is periodic of period $6$ and has a twisting
given by a class in $H^4(X,\bZ_{(2)})$, which the authors of
\cite{MR2681765,MR2810946,MR3431663} relate
to the notion of ``twisted String structures,''  which bear the same
relation to M-theory that that twisted spin$^c$-structures do in relation
to string theory.  (Just as anomaly cancellation for D-branes in the
presence of an $H$-flux requires D-branes to take their charges in
twisted $K$-theory, a similar anomaly cancellation in M-theory
suggests that branes should have charges in twisted $K(2)$.)

Another example of a twisted cohomology theory 
that comes up in M-theory is twisted cohomotopy ---
see \cite{MR4121615,MR4242295,MR4252880,MR4318438}.  Note that stable
cohomotopy is in some sense the universal cohomology theory.  (Unstable)
cohomotopy is not exactly a cohomology theory; in fact cohomotopy sets
are not always groups, but it reduces to stable cohomotopy in some
circumstances.  And calculations suggest that anomaly cancellation in
M-theory should force charges to live in twisted cohomotopy groups.

\bibliographystyle{amsplain}
\bibliography{twistedcohom}
\end{document}